\documentclass{article}

\usepackage{amsmath,amsfonts,amssymb,a4,enumerate,epsfig,theorem,psfrag}
\usepackage[latin1]{inputenc}

\newcommand{\Upper}{{\textrm{\rm Upper}}}

\newcommand{\RR}{{\mathbb{R}}}

\newcommand{\Ii}{{\cal{I}}}
\newcommand{\Jj}{{\cal{J}}}

\newcommand{\Pp}{{\cal{P}}}

\newcommand{\Ss}{{\mathbb{S}}}

\newcommand{\Uu}{{\cal{U}}}

\newcommand{\ILa}{\Ii_\Lambda}

\newcommand{\JLa}{\Jj_\Lambda}
\newcommand{\JLao}{\Jj^0_\Lambda}

\newcommand{\PLa}{\Pp_\Lambda}
\newcommand{\UpLa}{\Uu^\pi_\Lambda}

\newcommand{\llangle}{\langle\!\langle}
\newcommand{\rrangle}{\rangle\!\rangle}

\newcommand{\proof}{\smallskip{\noindent\bf Proof: }}

\def\qed{\unskip\nobreak\hfil\penalty50\hskip1.75em\null\nobreak\hfil
$\blacksquare$ {\parfillskip=0pt \finalhyphendemerits=0 \par}\medbreak}

\newcommand\capsize{\relax}

\newcommand\diag{\operatorname{diag}}

\newtheorem{lemma}{Lemma}[section]

\newtheorem{prop}[lemma]{Proposition}

\title{Reconstruction of tridiagonal matrices \\ from spectral data}
\author{Ricardo S. Leite, Nicolau C. Saldanha and Carlos Tomei
\footnote{The authors acknowledge support from CNPq, IM-AGIMB and FAPERJ.}
}

\begin{document}
\maketitle

\begin{abstract}
Jacobi matrices are parametrized by their eigenvalues and norming constants
(first coordinates of normalized eigenvectors):
this coordinate system breaks down at reducible tridiagonal matrices.
The set of real symmetric tridiagonal matrices with prescribed
simple spectrum is a compact manifold, admitting an open covering
by open dense sets $\UpLa$ centered at diagonal matrices $\Lambda^\pi$,
where $\pi$ spans the permutations.
{\it Bidiagonal coordinates} are a variant of norming constants
which parametrize each open set $\UpLa$ by the Euclidean space.

The reconstruction of a Jacobi matrix from inverse data
is usually performed by an algorithm introduced by de Boor and Golub.
In this paper we present a reconstruction procedure
from bidiagonal coordinates and show how to employ it
as an alternative to the de Boor-Golub algorithm.
The inverse bidiagonal algorithm rates well 
in terms of speed and accuracy.
\end{abstract}

\medbreak

{\noindent\bf Keywords:} Jacobi matrix, inverse eigenvalue problem,
bidiagonal coordinates. 

\smallbreak

{\noindent\bf MSC-class:} 65F18; 15A29.

\section{Introduction}

Recall that a real tridiagonal symmetric matrix
\[ T = \begin{pmatrix} a_1 & b_1 & & & \\
b_1 & a_2 & b_2 & & \\ & b_2 & a_3 & \ddots & \\
& & \ddots & \ddots & b_{n-1} \\ & & & b_{n-1} & a_n \end{pmatrix} \]
is a {\it Jacobi matrix} if $b_i > 0$ for all $i$.
Jacobi matrices have simple spectrum and their eigenvectors 
have nonzero first and last coordinates.
Thus, a Jacobi matrix $T$ diagonalizes uniquely as $T = Q^\ast \Lambda Q$,
$\Lambda = \diag(\lambda_1 < \cdots < \lambda_n)$,
provided we demand that the {\it norming constants}
$w_i = Q_{i1}$ be positive for all $i$.
Let $\JLao$ be the set of Jacobi matrices with given simple spectrum $\Lambda$
and $\Ss^{n-1}_+ =
\{(w_1, \ldots, w_n) \;|\; w_i > 0, w_1^2 + \cdots + w_n^2 = 1 \} $,
the open positive octant of the unit sphere $\Ss^{n-1} \subset \RR^n$.
Define the {\it map of norming constants}
$\omega_\Lambda: \JLao \to \Ss^{n-1}_+$
by $\omega_\Lambda(T) = w = (w_1, \ldots, w_n)$.
Moser (\cite{Moser}) proved that $\omega_\Lambda$ is a diffeomorphism.

On another route, numerical analysts considered the problem
of reconstructing a Jacobi matrix $T$ from its eigenvalues $\lambda_i$,
$i = 1\ldots n$, and the eigenvalues $\mu_i$, $i = 1\ldots n-1$,
of its bottom principal $(n-1) \times (n-1)$ minor.
The interlacing theorem requires that
$\lambda_1 < \mu_1 < \lambda_2 < \mu_2 < \cdots < \mu_{n-1} < \lambda_n$.
Existence and uniqueness of $T$, continuous dependence on
$\lambda$'s and $\mu$'s and an iterative algorithm to obtain $T$
were obtained by Hochstadt, Gray, Wilson and Hald
(\cite{Hochstadt}, \cite{Hochstadt2}, \cite{GW}, \cite{Hald});
A direct, stable algorithm was obtained by de Boor and Golub (\cite{BG}).
The first step of their algorithm is the computation
of the norming constants $w$ in terms of $\lambda$'s and $\mu$'s.
The problem then boils down to computing the inverse map
$\omega_\Lambda^{-1}(w)$,
as will be described in section 2.
For a survey of the Jacobi reconstruction problem, see \cite{EP}.

Let $\JLa$ be the closure of $\JLao$, the set of tridiagonal
symmetric matrices $T$ with spectrum $\Lambda$ and $b_i \ge 0$.
The map of norming constants $\omega_\Lambda$ extends continuously to $\JLa$
but this extension is no longer injective.
Indeed, for $n = 3$ and $\Lambda = \diag(1,2,4)$ we have
\[ \begin{pmatrix} 1 & 0 & 0 \\ 0 & 3 - \cos 2t & \sin 2t \\
0 & \sin 2t & 3 + \cos 2t \end{pmatrix} =
\begin{pmatrix} 1 & 0 & 0 \\ 0 & \cos t & \sin t \\
0 & -\sin t & \cos t \end{pmatrix}
\Lambda
\begin{pmatrix} 1 & 0 & 0 \\ 0 & \cos t & -\sin t \\
0 & \sin t & \cos t \end{pmatrix} \]
and therefore $\omega_\Lambda(T) = (1,0,0)$, $\mu_1 = 4$ and $\mu_2 = 6$ 
for all such $T$.
Thus, in some sense, any reconstruction algorithm either
from $\lambda$'s and $\mu$'s or from $\lambda$'s and $w$'s
must degenerate at some points of the boundary of $\JLa$.

In \cite{LST}, the authors introduced {\it bidiagonal coordinates},
a variant of norming constants which behaves well at the boundary.
In this paper we provide a direct reconstruction algorithm
from bidiagonal coordinates with good behavior at boundary points,
where norming constants break down.
The conversion from norming constants to bidiagonal coordinates is
simple, and the resulting algorithm is comparable in time and space
with that of de Boor and Golub, being more accurate in many cases.

\begin{figure}[ht]
\begin{center}
\epsfig{height=50mm,file=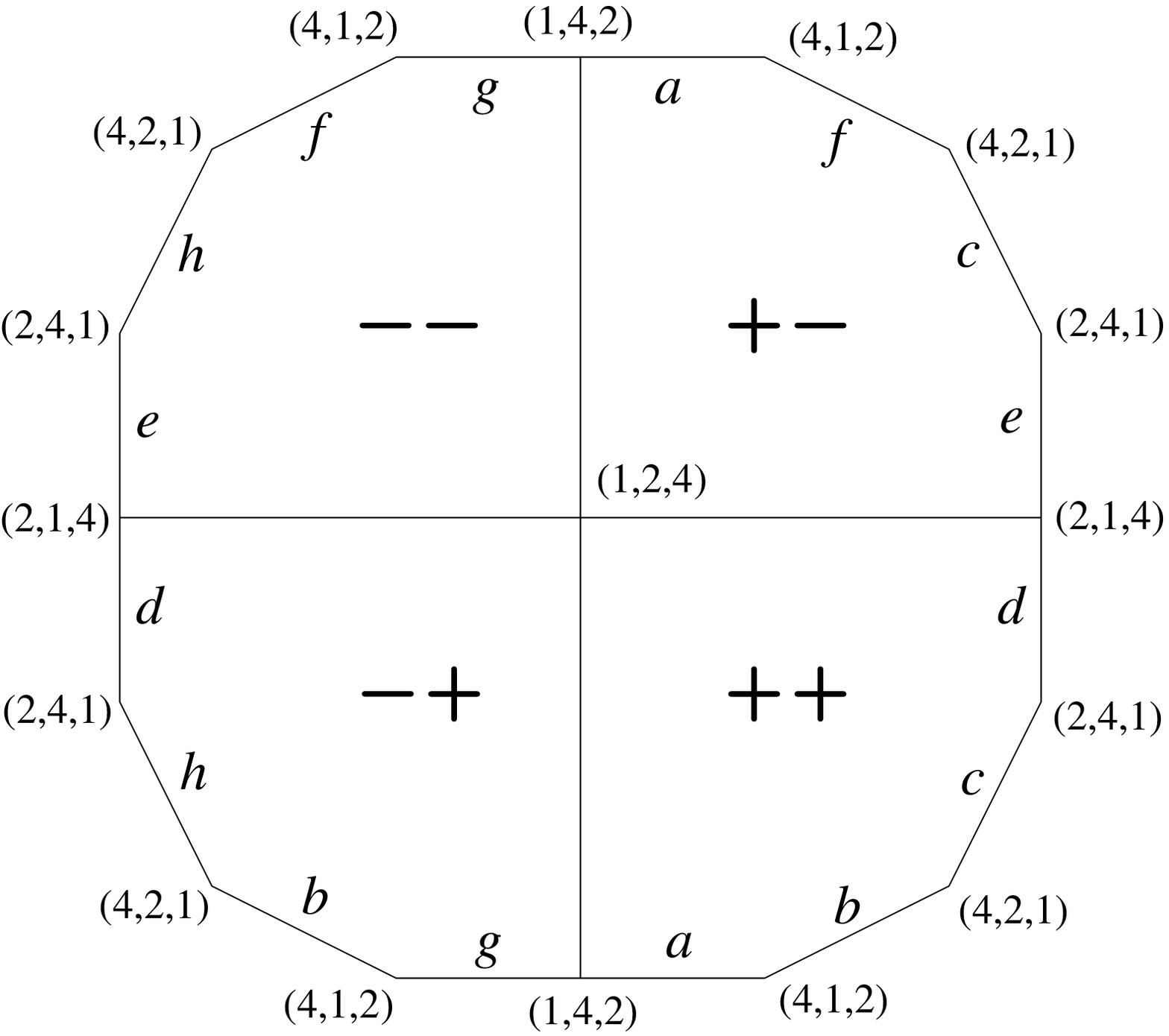}
\end{center}
\caption{\capsize The manifold $\ILa$ for $\Lambda = \diag(1,2,4)$}
\label{fig:ilailah}
\end{figure}

Matters become clearer with some geometric vocabulary.
Let $\ILa \supset \JLao$ be the set of tridiagonal symmetric matrices
with spectrum $\Lambda$. As proved in \cite{Tomei},
$\ILa$ is a compact oriented manifold of dimension $n-1$.
The closure $\JLa \subset \ILa$ of $\JLao$ is homeomorphic 
to the convex polytope $\PLa$ with vertices
$\Lambda^\pi = \diag(\lambda_{\pi(1)}, \ldots, \lambda_{\pi(n)})$
where $\pi$ spans the set of permutations of $\{1, 2, \ldots, n\}$
(\cite{Tomei}, \cite{BFR}).
Each of the $2^{n-1}$ possible choices of signs for the entries
$b_i$ define a closed subset of $\ILa$ which is isomorphic to $\JLa$:
as is well known, dropping the signs of the off-diagonal entries
of a tridiagonal symmetric matrix does not change its spectrum.
Thus, $\ILa$ can be constructed by glueing $2^{n-1}$ copies of $\PLa$
along faces consisting of reducible tridiagonal matrices
(i.e., matrices for which some $b_i$ is zero).
In figure \ref{fig:ilailah} we show what happens for $n = 3$.
The polytope $\PLa$ is a hexagon and
in each of its copies we indicate the signs of $b_1$ and $b_2$.
Vertices are diagonal matrices and edges consist of reducible matrices;
edges with the same label are glued.
It follows that $\ILa$ is a bitorus for $n=3$.

Removal of the outer boundary edges in the picture yields an open dense
subset of $\ILa$ centered at $\Lambda$.
As we shall see, bidiagonal coordinates can be smoothly defined
on this set, yielding an explicit diffeomorphism with $\RR^2$.
More generally, for each $\Lambda$ and each permutation $\pi$, 
we define an open dense subset $\UpLa$ of $\ILa$.
The complement of $\UpLa$ consists of matrices $T$ for which
there exist $i \le k < n$ with $b_k = 0$ and
$\lambda_{\pi(i)}$ belongs to the spectrum of the bottom
principal $(n-k) \times (n-k)$ minor.
As in the example above, $\UpLa$ is centered at $\Lambda^\pi$
and bidiagonal coordinates provide a diffeomorphism
between $\UpLa$ and $\RR^{n-1}$.
Also, the sets $\UpLa$ form an open cover of $\ILa$:
this is the crucial property for the local study
of iterations preserving tridiagonality and spectrum
as performed in \cite{LST}.


\section{The de Boor-Golub algorithm}

In this section we present part of the contents of \cite{BG}
phrased in such a way as to emphasize the differences and similarities
between this more well known algorithm and the inverse bidiagonal algorithm,
to be presented in section 4.
We assume that the off-diagonal entries $b_i$
of $T$ are positive so that $T$ is a Jacobi matrix.
Write $T = Q^\ast \Lambda Q$ where 
$\Lambda = \diag(\lambda_1, \ldots, \lambda_n)$
(in this section, the simple eigenvalues $\lambda_i$
are taken in an arbitrary order),
and $Q$ is orthogonal with positive first column.
For $D_b = \diag(1, b_1, b_1b_2, \ldots, b_1b_2\cdots b_{n-1})$
and $D_w = \diag(Q_{11}, Q_{21}, \ldots, Q_{n1})
= \diag(w_1, w_2, \ldots, w_n)$, write
$\tilde P = D_b Q^\ast D_w^{-1}$ so that
\[ \tilde T = D_b T D_b^{-1} = \begin{pmatrix} a_1 & 1 & & & \\
b_1^2 & a_2 & 1 & & \\ & b_2^2 & a_3 & \ddots & \\
& & \ddots & \ddots & 1 \\ & & & b_{n-1}^2 & a_n \end{pmatrix} =
\tilde P \Lambda {\tilde P}^{-1}. \]
Let $\tilde p^\ast_{k-1} = e^\ast_k \tilde P$
be the $k$-th row of $\tilde P$;
in particular, $\tilde p^\ast_0 = (1, 1, \ldots, 1)$.
We have $\tilde T \tilde P = \tilde P \Lambda$ and the rows
of $\tilde P$ satisfy the recursion
\begin{align}
\tilde p^\ast_1 &= \tilde p^\ast_0 \Lambda - a_1 \tilde p^\ast_0, \notag\\
\tilde p^\ast_{k+1} &= \tilde p^\ast_k \Lambda -
a_{k+1} \tilde p^\ast_k - b_k^2 \tilde p^\ast_{k-1}, \quad
0 < k < n-1, \notag\\
0 &= \tilde p^\ast_{n-1} \Lambda -
a_n \tilde p^\ast_{n-1} - b_{n-1}^2 \tilde p^\ast_{n-2}. \notag
\end{align}
Furthermore, the vectors $\tilde p_k$ form an orthogonal basis
under the inner product
$\llangle u, v \rrangle = \langle u, D_w^2 v \rangle$.
For $0 \le k < n$, let $\check p_k$ be the unique polynomial
of degree less than $n$ satisfying $\check p_k(\lambda_j) = (\tilde p_k)_j$:
we have $\check p_0 = 1$,
\[ \check p_1 = t \check p_0 - a_1 \check p_0, \quad
\check p_{k+1} = t \check p_k - a_{k+1} \check p_k - b_k^2 \check p_{k-1},
\quad 0 < k < n, \]
so that $\check p_k$ is a monic polynomial of degree $k$.

In the notation of \cite{BG}, let $p_k(t) = \det(tI - T_k)$ where
$T_k$ is the principal minor of $T$ consisting of the first $k$
rows and columns; set also $p_0 = 1$.
The expansion of the determinant along the last row of each minor yields
\begin{equation}
p_1 = t p_0 - a_1 p_0, \quad
p_{k+1} = t p_k - a_{k+1} p_k - b_k^2 p_{k-1}, \quad 0 < k < n
\label{eq:recp} \end{equation}
and therefore $p_k = \check p_k$ since $p_0 = \check p_0$
and both sequences satisfy the same recurrence.
Equivalently, the $j$-th coordinate of the vector $\tilde p^\ast_k$
is $p_k(\lambda_j)$.

Summing up, assume $\Lambda$ and $D_w = \diag(w_1, \ldots, w_n)$ given. 
The linear bijection between the space of real polynomials
of degree less than $n$ and $\RR^n$ given by evaluation on the $\lambda_j$'s
allows us to pull back the inner product $\llangle\cdot,\cdot\rrangle$ 
giving rise to an inner product on polynomials:
\[ \llangle q_1, q_2 \rrangle =
\sum_{j=1}^n w_j^2 q_1(\lambda_j) q_2(\lambda_j). \]
The de Boor-Golub algorithm now constructs the monic polynomials $p_k$
using the orthogonality condition $\llangle p_k, p_{k'} \rrangle = 0$
for $k \ne k'$.
Recursion \ref{eq:recp} for the polynomials $p_k$ obtains
the entries of $T$.


\section{Bidiagonal coordinates}

We quote some of the results and notations of \cite{LST}
to be used in the inverse bidiagonal algorithm.
Diagonalize $T \in \ILa$ as $T = Q^\ast \Lambda Q$,
$\Lambda = \diag(\lambda_1 < \lambda_2 < \cdots < \lambda_n)$,
and factor $Q = PLU$ where $P$ is a permutation matrix,
$L$ is lower unipotent and $U$ is upper triangular.
For a permutation $\pi$, let $P_\pi$ be the permutation matrix
with $(i,j)$ entry equal $1$ iff $i = \pi(j)$
(thus $P_{\pi_1 \pi_2} = P_{\pi_1} P_{\pi_2}$ and $P_\pi e_i = e_{\pi(i)}$).
The $PLU$ factorization can usually be done
for several choices of the permutation matrix:
it turns out that we can take $P = P_\pi$ if and only if $T \in \UpLa$,
the open dense subset of $\ILa$ presented in the introduction.
Write $Q_\pi = E P_\pi^{-1} Q = L_\pi U_\pi$ where $E$ is a diagonal
matrix with diagonal entries equal to $1$ and $-1$,
$L_\pi$ is lower unipotent and $U_\pi$ is upper triangular
with positive diagonal.
The rows of $Q_\pi$ are eigenvectors of $T$ but
their first coordinates are not necessarily nonnegative:
instead, signs are determined from the fact that
the determinants of leading principal minors of $Q_\pi$ are positive.

Let $B_\pi = L_\pi^{-1} \Lambda L_\pi = R_\pi^{-1} T R_\pi$
where $R_\pi = U_\pi^{-1}$ so that $L_\pi = Q_\pi R_\pi$.
From the first formula, $B_\pi$ is lower triangular;
from the second, it is upper Hessenberg;
thus, $B_\pi$ is lower bidiagonal:
\[ B_\pi = \begin{pmatrix} \lambda^\pi_1 & & & & \\
\beta^\pi_1 & \lambda^\pi_2 & & & \\ & \beta^\pi_2 & \lambda^\pi_3 & & \\
& & \ddots & \ddots & \\ & & & \beta^\pi_{n-1} & \lambda^\pi_n \end{pmatrix}. \]
The map $\psi_\pi: \UpLa \to \RR^{n-1}$ taking $T$ to
the {\it $\pi$-bidiagonal coordinates}
$(\beta^\pi_1, \ldots, \beta^\pi_{n-1})$ is a diffeomorphism.
Indeed, start from an explicit formula for the matrix $L_\pi$
in terms of bidiagonal coordinates:
\[ L_\pi = \begin{pmatrix}
1 & 0 & 0 & \cdots & 0 \\ \\
\frac{\beta^\pi_1}{\lambda^\pi_2 - \lambda^\pi_1} & 1 & 0 & \cdots & 0 \\ \\
\frac{\beta^\pi_1 \beta^\pi_2}%
{(\lambda^\pi_3 - \lambda^\pi_1)(\lambda^\pi_3 - \lambda^\pi_2)} &
\frac{\beta^\pi_2}{\lambda^\pi_3 - \lambda^\pi_2}  & 1 & & 0 \\
\vdots & \vdots & & \ddots & \\
\frac{\beta^\pi_1 \beta^\pi_2 \cdots \beta^\pi_{n-1}}%
{(\lambda^\pi_n - \lambda^\pi_1)(\lambda^\pi_n - \lambda^\pi_2)\cdots%
(\lambda^\pi_n - \lambda^\pi_{n-1})} &
\frac{\beta^\pi_2 \cdots \beta^\pi_{n-1}}%
{(\lambda^\pi_n - \lambda^\pi_2)\cdots (\lambda^\pi_n - \lambda^\pi_{n-1})} &
& & 1
\end{pmatrix}.
\]
Given $L_\pi$, its $QR$ factorization yields $Q_\pi$ and $R_\pi$,
from which one obtains $T = R_\pi B_\pi R_\pi^{-1}$.
A straightforward computation shows that $b_i$ and $\beta^\pi_i$
have the same sign and that near a diagonal matrix,
$\beta^\pi_i$ equals $b_i$ to first order;
more, for the inverse map
$\phi_\pi = (\psi_\pi)^{-1}: \RR^{n-1} \to \UpLa \subset \ILa$,
\[ \phi_\pi(0, \ldots, 0) + D\phi_\pi(0, \ldots, 0) (u_1, \ldots, u_{n-1}) =
\begin{pmatrix} \lambda^\pi_1 & u_1 & & \\
u_1 & \lambda^\pi_2 & u_2 & \\ & u_2 & \lambda^\pi_3 & \\ & & & \ddots
\end{pmatrix}. \]

For any permutation $\pi$ and any $T \in \UpLa$,
the norming constants $w^\pi_i = w_{\pi(i)}$ and the
$\pi$-bidiagonal coordinates $\beta^\pi_i$ are related by 
\[ w^\pi_i = w^\pi_1 \frac{\beta^\pi_1 \cdots \beta^\pi_{i-1}}%
{(\lambda^\pi_i - \lambda^\pi_1) \cdots (\lambda^\pi_i - \lambda^\pi_{i-1})},
\quad 2 \le i \le n, \]
\[ \beta^\pi_i = \frac{(\lambda^\pi_{i+1} - \lambda^\pi_1)\cdots
(\lambda^\pi_{i+1} - \lambda^\pi_i) w^\pi_{i+1}}
{(\lambda^\pi_{i} - \lambda^\pi_1)\cdots
(\lambda^\pi_{i} - \lambda^\pi_{i-1}) w^\pi_{i}},
\quad 1 \le i \le n-1. \]

The inverse bidiagonal algorithm, presented in the next section,
obtains the matrix $R_\pi$ in another way,
closer in spirit to the recursions in section 2.
The basic version of this algorithm receives as input a permutation $\pi$,
eigenvalues $\lambda^\pi_i$ and bidiagonal coordinates $\beta^\pi_i$
and returns the corresponding tridiagonal matrix $T \in \UpLa$.
The de Boor-Golub algorithm, instead, receives as input
the eigenvalues $\lambda_i$ and the norming constants $w_i$:
in this case a simultaneous permutation $\pi$ is innocuous,
at least with exact arithmetic.

\section{The inverse bidiagonal algorithm}

We first describe a preliminary version of the algorithm,
which only works in the irreducible case,
where all $\beta^\pi_k$ (or, equivalently, all $b_k$) are nonzero.
Write $\tilde T = \hat R B_\pi \hat R^{-1}$ where $\hat R = D_b R_\pi$
is an upper triangular matrix with rows $\hat r_k^\ast$.
Clearly, $\hat r_1^\ast = e_1^\ast \hat R = e_1^\ast D_b R_\pi =
e_1^\ast R_\pi =
e_1^\ast Q_\pi^\ast L_\pi = (Q_\pi e_1)^\ast L_\pi =
(L_\pi U_\pi e_1)^\ast L_\pi =
u_{11} (L_\pi e_1)^\ast L_\pi = u_{11} e_1^\ast L_\pi^\ast L_\pi$
and therefore
$\hat r_1 = c L_\pi^\ast L_\pi e_1$, the value of $c = u_{11} > 0$
being irrelevant throughout the algorithm.

Equate rows in $\tilde T \hat R = \hat R B_\pi$,
\[ \begin{pmatrix} a_1 & 1 & & & \\
b_1^2 & a_2 & 1 & & \\ & b_2^2 & a_3 & \ddots & \\
& & \ddots & \ddots & 1 \\ & & & b_{n-1}^2 & a_n \end{pmatrix} 
\begin{pmatrix} \hat r_1^\ast \\ \hat r_2^\ast \\ \hat r_3^\ast \\ \vdots
\\ \hat r_n^\ast \end{pmatrix} = 
\begin{pmatrix} \hat r_1^\ast \\ \hat r_2^\ast \\ \hat r_3^\ast \\ \vdots
\\ \hat r_n^\ast \end{pmatrix} 
\begin{pmatrix} \lambda^\pi_1 & & & & \\
\beta^\pi_1 & \lambda^\pi_2 & & & \\ & \beta^\pi_2 & \lambda^\pi_3 & & \\
& & \ddots & \ddots & \\
& & & \beta^\pi_{n-1} & \lambda^\pi_n \end{pmatrix}, \]
to obtain $\hat r_{k+1}^\ast = \hat r_k^\ast B_\pi
- a_k \hat r_k^\ast - b_{k-1}^2 \hat r_{k-1}^\ast$.
Since $\hat R$ is known to be upper triangular,
this recursion, together with the initial term $\hat r^\ast_1$,
allows us to compute the coefficients $a_i$, $i = 1, \ldots, n$
and $b^2_i$, $i = 1, \ldots, n-1$.
The numbers $b_i$ and $\beta^\pi_i$ have the same sign in $\UpLa$:
this completes the preliminary version of the reconstruction algorithm
for irreducible matrices.


We need to modify the algorithm in order to extend it to the general case.
For an integer $k \ge 0$, set
\begin{equation}
L_{\pi,k} = \begin{pmatrix}
1 & 0 & 0 & \cdots & 0 \\ \\
\frac{(\beta^\pi_1)^k}{\lambda^\pi_2 - \lambda^\pi_1} & 1 & 0 & \cdots & 0 \\ \\
\frac{(\beta^\pi_1 \beta^\pi_2)^k}%
{(\lambda^\pi_3 - \lambda^\pi_1)(\lambda^\pi_3 - \lambda^\pi_2)} &
\frac{(\beta^\pi_2)^k}{\lambda^\pi_3 - \lambda^\pi_2}  & 1 &  & 0 \\
\vdots & \vdots & & \ddots & \\
\frac{(\beta^\pi_1 \beta^\pi_2 \cdots \beta^\pi_{n-1})^k}%
{(\lambda^\pi_n - \lambda^\pi_1)(\lambda^\pi_n - \lambda^\pi_2)\cdots%
(\lambda^\pi_n - \lambda^\pi_{n-1})} &
\frac{(\beta^\pi_2 \cdots \beta^\pi_{n-1})^k}%
{(\lambda^\pi_n - \lambda^\pi_2)\cdots (\lambda^\pi_n - \lambda^\pi_{n-1})} &
& & 1
\end{pmatrix}
\label{eq:Lk}
\end{equation}
and $B_{\pi,k} = L_{\pi,k}^{-1} \Lambda^\pi L_{\pi,k}$, or, more explicitly,
\[ B_{\pi,k} =
\begin{pmatrix} \lambda^\pi_1 & & & & \\
(\beta^\pi_1)^k & \lambda^\pi_2 & & & \\
& (\beta^\pi_2)^k & \lambda^\pi_3 & & \\
& & \ddots & \ddots & \\ & & & (\beta^\pi_{n-1})^k & \lambda^\pi_n
\end{pmatrix}. \]
Still in the irreducible case,
define $D_\beta = \diag(1,\beta^\pi_1, \beta^\pi_1\beta^\pi_2, \ldots,
\beta^\pi_1\beta^\pi_2 \cdots \beta^\pi_{n-1})$ and
$\tilde R = c^{-1} \hat R D_\beta^{-1}$ with rows $\tilde r^\ast_k$
so that $B_{\pi,2} = D_\beta B_\pi D_\beta^{-1}$ and
$\tilde T \tilde R = \tilde R B_{\pi,2}$.
Straightforward computations verify that
$\tilde r_1 = L_{\pi,2}^\ast L_{\pi,0} e_1$.
Expanding the matrix products as above we obtain the recursion
$\tilde r_{k+1}^\ast = \tilde r_k^\ast B_{\pi,2}
- a_k \tilde r_k^\ast - b_{k-1}^2 \tilde r_{k-1}^\ast$.
Thus, from $\tilde r_{k-1}^\ast$ and $\tilde r_k^\ast$
we compute $\tilde r_k^\ast B_{\pi,2}$,
then $b_{k-1}$ and $a_k$ and finally $\tilde r_{k+1}^\ast$.
This completes the description of the {\it inverse bidiagonal algorithm}
for irreducible matrices;
we now prove that this procedure works for any $\beta^\pi \in \RR^{n-1}$,
obtaining all matrices $T \in \UpLa$.
Let $\Upper^+(\RR,n)$ be
the group of upper triangular matrices with positive diagonal.

\begin{prop}
\label{prop:rho}
There is a smooth function $\rho: \UpLa \to \Upper^+(\RR,n)$
satisfying $\rho(T) = (R_\pi)_{11} D_b R_\pi D_\beta^{-1}$
for all irreducible matrices $T \in \UpLa$.
\end{prop}

Here, as in section 3, $T = Q_\pi^\ast \Lambda^\pi Q_\pi$,
$L_\pi = Q_\pi R_\pi$,
$Q_\pi$ orthogonal, $L_\pi$ lower unipotent
and $R_\pi \in \Upper^+(\RR,n)$.
The purpose of this proposition is to make sense of $\tilde R$
for reducible matrices $T$ (or, equivalently, for $\beta^\pi$
with some zero coordinate).
The formula in the statement defines $\rho(T)$ as $\tilde R$
for irreducible $T$ but otherwise involves divisions by zero.


\proof
Define $\tilde\rho: \UpLa \to \RR^{n\times n}$ row by row:
let $\tilde\rho_k^\ast$ denote the $k$-th row of $\tilde\rho(T)$
and set $\tilde\rho_1 = L_{\pi,2}^\ast L_{\pi,0} e_1$,
$\tilde\rho_{k+1}^\ast = \tilde\rho_k^\ast B_{\pi,2}
- a_k \tilde\rho_k^\ast - b_{k-1}^2 \tilde\rho_{k-1}^\ast$.
The function $\tilde\rho$ is clearly smooth in $\UpLa$.
Also, we proved above that $\tilde\rho(T) = \rho(T) = \tilde R$
for irreducible $T$.
Thus, by continuity, $\tilde\rho(T)$ is always upper triangular
with nonnegative diagonal entries $\tilde\rho_{k,k}$.
In the irreducible case,
\[ \tilde\rho_{k,k} = 
\frac{b_1 b_2 \cdots b_{k-1}}%
{c \beta^\pi_1 \beta^\pi_2 \cdots \beta^\pi_{k-1}} (R_\pi)_{k,k}. \]
It remains to prove that $\tilde\rho_{k,k} \ne 0$ for reducible $T$
so that we can then set $\rho = \tilde\rho$.

One way of completing the proof is recalling from \cite{LST}
that the quotients $b_j/\beta^\pi_j$ are smooth positive functions on $\UpLa$.
Alternatively, from $B_{\pi,2} = L_{\pi,2}^{-1} \Lambda L_{\pi,2}$,
we can write $(B_{\pi,2}^\ast)^{k-1} \tilde\rho_1 =
L_{\pi,2}^\ast \Lambda^{k-1} L_{\pi,0} e_1$.
The coordinates of $L_{\pi,0} e_1$ are all nonzero by equation \ref{eq:Lk}
and, from the standard Vandermonde argument,
the vectors $\Lambda^{k-1} L_{\pi,0} e_1$, $k = 1, \ldots, n$, form a basis;
since $L_{\pi,2}^\ast$ is invertible,
so do the vectors $(B_{\pi,2}^\ast)^{k-1} \tilde\rho_1$.
From the recursion formula,
so do the vectors $\tilde\rho_k$ and we are done.
\qed

In general, we start from $\tilde r_1^\ast$ and use the recursive formula
$\tilde r_{k+1}^\ast = \tilde r_k^\ast B_{\pi,2}
- a_k \tilde r_k^\ast - b_{k-1}^2 \tilde r_{k-1}^\ast$.
More precisely, assume by induction that
$\tilde r_{k-1}^\ast$ and $\tilde r_k^\ast$ are known.
From the proposition, $\tilde r_{k,k}$ and $\tilde r_{k-1,k-1}$
are positive. The first nonzero coordinate of $\tilde r_k^\ast B_{\pi,2}$
occupies position $k-1$ and equals $(\beta^\pi_{k-1})^2 \tilde r_{k,k}$.
The algorithm then calculates
$b_{k-1} = \beta^\pi_{k-1} \sqrt{\tilde r_{k,k}/\tilde r_{k-1,k-1}}$,
so that the square root is evaluated at a strictly positive number.
Notice that the algorithm treats uniformly all $\beta^\pi \in \RR^{n-1}$,
i.e., there is no checking of signs or division into cases.
The values of $a_k$ and of $\tilde r_{k+1}^\ast$,
the $(k+1)$-th row of the upper triangular matrix $\tilde R$,
are now easily obtained,
concluding the computation of $\phi_\pi(\beta^\pi_1, \ldots, \beta^\pi_{n-1})$
and the description of the inverse bidiagonal algorithm.


\section{Accuracy and tight permutations}

Empirical evidence indicates that a good choice of the permutation $\pi$
is extremely important for the accuracy of the inverse bidiagonal algorithm.
One is reminded of Gaussian elimination,
where pivoting strategies have a similar effect.
There is a crucial difference, however.
In Gaussian elimination, the permutation is chosen
along the process; the inverse bidiagonal algorithm admits
no easy way to accomodate a change of permutation in mid-flight.
As to estimating accuracy with respect to the choice of permutation,
our theoretical understanding is limited and we provide instead
a simple numerical experiment
\footnote{Maple worksheets for all experiments in this paper are
available at \hfil\break
{\tt http://www.mat.puc-rio.br/\~\relax nicolau/papers/invbi-mw}.}.
We start with random inverse data for an $8 \times 8$ matrix
and perform the inverse bidiagonal algorithm
for each of the $8!$ permutations with $8$ digits of precision.
Results are then compared with the ``correct'' answer,
computed with an exaggerated number of digits.
Different permutations yield very different errors:
the smallest error is $2.0 \cdot 10^{-7}$,
there are $19$ other permutations with error smaller than $3 \cdot 10^{-7}$
and there are $8$ permutations with error greater then $7 \cdot 10^{-2}$.
The error here is defined as
\[ \sum_{i = 1}^n |a_i - \tilde a_i| + \sum_{i=1}^{n-1} |b_i - \tilde b_i| \]
where $a_i$ and $b_i$ are the ``correct'' values and 
$\tilde a_i$ and $\tilde b_i$ are the computed values.
All entries of $T$ have absolute value smaller than $1$.
In this section, we present a strategy for choosing $\pi$.

Let $\tau_i$ be the transposition $(k, k+1)$ in cycle notation.
Two permutations $\pi_0$ and $\pi_1$ {\it differ by $\tau_k$}
if $\pi_1 = \pi_0 \circ \tau_k$.
Thus, $\Lambda^{\pi_1}$ is obtained from $\Lambda^{\pi_0}$
by interchanging the $(k,k)$ and $(k+1,k+1)$ entries.
Bidiagonal coordinates $\beta^{\pi_0}_i$ and $\beta^{\pi_1}_i$
are equal except for
\begin{equation}
\beta^{\pi_1}_{k-1} = q^{\pi_0}_k \beta^{\pi_0}_{k-1}, \quad
\beta^{\pi_1}_{k} = - (q^{\pi_0}_k)^{-2} \beta^{\pi_0}_{k}, \quad
\beta^{\pi_1}_{k+1} = - q^{\pi_0}_k \beta^{\pi_0}_{k+1}, 
\label{eq:newbeta} \end{equation}
where
\[ q^{\pi_0}_k =
\frac{\beta^{\pi_0}_k}{\lambda^{\pi_0}_{k+1} - \lambda^{\pi_0}_k}, \]
as can be proved from the formulae relating $\beta$'s and $w$'s
in section 3. 

Given inverse data $\lambda_i$ and $w_i$,
we call a permutation $\pi$ {\it tight} if $|q^\pi_i| \le 1$
for all $i = 1, \ldots, n-1$ and we say that the transposition
$\tau_k$ is {\it $\pi$-tightening} if $|q^\pi_k| > 1$.
From equation \ref{eq:newbeta}, it is easy to see that 
if $\tau_k$ is $\pi$-tightening then $|q^{\pi \circ \tau_k}_k| < 1$.
Clearly, $\pi$ is tight if and only if there are no
$\pi$-tightening transpositions.
A {\it tightening sequence} is a maximal sequence
$(\pi_m)$ of permutations such that $\pi_{m+1}$ differs
from $\pi_m$ by a $\pi_m$-tightening transposition $\tau_{k_m}$.
Thus, a tightening sequence is either infinite
or ends at a tight permutation.

\begin{lemma}
\label{lemma:tight}
Tightening sequences are finite.
\end{lemma}

\proof
Assume by contradiction that there exists an infinite tightening sequence:
clearly, there exist $m_0 < m_1$ with $\pi_{m_0} = \pi_{m_1}$.
We show that there are no such cycles.

Set $p_{k,\pi} = \prod_{i \ge k} |\beta^\pi_i|$ and
use the lexicographical order to define a total order
in the permutation group:
$\pi_0 \prec \pi_1$ if and only if there exists $k'$
such that $p_{k,\pi_0} = p_{k,\pi_1}$ for $k < k'$
and $p_{k',\pi_0} < p_{k',\pi_1}$.
For any $m$, it follows from equation \ref{eq:newbeta}
that $p_{k,\pi_{m+1}} = p_{k,\pi_m}$ for $k < k_m$
and $p_{k_m,\pi_{m+1}} < p_{k_m,\pi_m}$,
implying $\pi_{m+1} \prec \pi_m$.
By transitivity, $\pi_{m_0} \prec \pi_{m_1} = \pi_{m_0}$,
a contradiction.
\qed

In the example discussed above,
there were $4$ tight permutations with errors between
$2.8 \cdot 10^{-7}$ and $5.3 \cdot 10^{-7}$.
Empirical evidence shows that this is frequent:
tight permutations usually yield
small errors in the inverse bidiagonal algorithm.
Thus, upon receiving inverse data $\lambda_i$ and $w_i > 0$
for a Jacobi matrix $T$, we first order the $w$'s in decreasing order
to obtain a permutation $\pi_0$ and then apply tightening transpositions
until we reach a tight permutation.
Experiments suggest that this takes approximately $n/2$ sweeps.

\section{Operational costs}

We estimate the number of operations (or flops) and
the amount of memory necessary to execute the inverse bidiagonal algorithm.
As usual, we only keep track of the number of products and quotients.

The process of finding a tight permutation will not be carefully
examined: suffice it to say that, from empirical evidence,
the number of operations is approximately $C n^2$ where $C < 1/2$.

The matrices $B_{\pi,2}$ and $L_{\pi,2}$ will come up along the algorithm
and it is therefore convenient to compute and keep
the squares $(\beta_i^\pi)^2$,
with an initial cost of $n$ operations and $n$ storage units.
The first major step of the algorithm is
the computation of  $\tilde r_1 = L_2^\ast L_0 e_1$.
For $n = 4$, after reordering terms, $\tilde r_1$ becomes
\begin{align}
& \Bigg(
\frac{(\beta^\pi_1)^2(\beta^\pi_2)^2(\beta^\pi_3)^2}
{ (\lambda^\pi_4-\lambda^\pi_1)^2 (\lambda^\pi_4-\lambda^\pi_2)^2
( \lambda^\pi_4-\lambda^\pi_3)^2}
+ \frac{(\beta^\pi_1)^2(\beta^\pi_2)^2}
{ (\lambda^\pi_3-\lambda^\pi_1)^2 (\lambda^\pi_3-\lambda^\pi_2)^2}
+ \frac{(\beta^\pi_1)^2}{(\lambda^\pi_2-\lambda^\pi_1)^2}
+ 1,
\notag\\
&\phantom{ X}
\frac{(\beta^\pi_2)^2(\beta^\pi_3)^2}
{ (\lambda^\pi_4-\lambda^\pi_1) (\lambda^\pi_4-\lambda^\pi_2)^2
(\lambda^\pi_4-\lambda^\pi_3)^2 }
+ \frac{(\beta^\pi_2)^2}{ (\lambda^\pi_3-\lambda^\pi_1)
(\lambda^\pi_3-\lambda^\pi_2)^2 }
+ \frac{1}{(\lambda^\pi_2-\lambda^\pi_1)},
\notag\\
&\phantom{ X}
\frac {(\beta^\pi_3)^2}
{ (\lambda^\pi_4-\lambda^\pi_1) (\lambda^\pi_4-\lambda^\pi_2)
(\lambda^\pi_4-\lambda^\pi_3)^2 }
+ \frac{1}{ (\lambda^\pi_3-\lambda^\pi_1)  (\lambda^\pi_3-\lambda^\pi_2 ) },
\notag\\
&\phantom{ X}
\frac {1}
{ (\lambda^\pi_4-\lambda^\pi_1) (\lambda^\pi_4-\lambda^\pi_2)
(\lambda^\pi_4-\lambda^\pi_3) }
\Bigg). \notag
\end{align}
These terms are computed from bottom to top of each column,
following the obvious patterns,
with a cost of approximately $3n^2/2$ operations
and $n$ storage units.

The recursion formula which obtains $b_{k-1}$, $a_k$
and $\tilde r_{k+1}$ only requires $\tilde r_{k-1}$ and $\tilde r_k$
so that we only need to keep at most three rows of the triangular
matrix $\tilde R$ at any given time.
The number of operations is approximately $2n^2$;
also, $n-1$ square roots are needed.

Summing up, given a tight permutation $\pi$ and the values
of $\beta^\pi_i$, a run takes approximately $7n^2/2$ operations,
$n$ square roots and $4n$ storage units
(provided some units do double duty,
first as entries of $\tilde R$ and later as $a$'s or $b$'s).

\section{Benchmarks}

In this section, we compare the
de Boor-Golub and inverse bidiagonal algorithms in a few scenarios.
This is only possible for irreducible matrices
since otherwise, as we saw, the norming constants $w$ break down.
The inverse bidiagonal algorithm receives as input
permuted eigenvalues $\lambda^\pi_i$, $i = 1,\ldots, n$,
and bidiagonal coordinates $\beta^\pi_i$, $i = 1,\ldots, n-1$.
In order to allow for comparisons, we must step back and
provide as input the eigenvalues $\lambda_i$ and the norming constants $w_i$:
we then obtain a tight permutation $\pi$ and compute $\beta^\pi_i$.

It is a common feature of both algorithms that the coefficients
$a_i$, $i = 1, \ldots, n$ and $b_i$, $i = 1, \ldots, n-1$,
are obtained in the order $a_1, b_1, a_2, b_2, a_3, \ldots$.
Also, both algorithms admit a reversal by conjugation.
More precisely, let $P_\rho$ be the permutation matrix
with $(P_\rho)_{ij} = 1$ if and only if $i+j=n+1$.
Let $\lambda_i$ and $w_i$ be the inverse data for a Jacobi matrix $T$:
the inverse data for $\tilde T = P_\rho T P_\rho$ is $\lambda_i$ and
\[ \tilde w_i = \frac{c}{w_i \prod_{j \ne i} |\lambda_i - \lambda_j|} \]
for some positive normalizing constant $c$ (\cite{BG}).
From data $\lambda_i$ and $\tilde w_i$ either algorithm obtains,
in this order,
$\tilde a_1 = a_n, \tilde b_1 = b_{n-1}, \tilde a_2 = a_{n-1}, \ldots$.
Experiments show that it is far wiser to do both things,
i.e., to compute the top half of $T$ directly from $w_i$
and the bottom half from $\tilde w_i$.
In the examples below, this strategy, the {\it two-sided} algorithms,
is always adopted.

We implement in a Maple worksheet both the two-sided de Boor-Golub (BG) and
the two-sided inverse bidiagonal (BI) algorithms, generate a sequence
of random inverse problems and compare errors 
for different values of the dimension and of
the number of significative digits.

In the first class of examples,
random real symmetric tridiagonal matrices $T$ are obtained as follows:
the nonzero entries are independent random variables
with a Gaussian distribution centered at $0$ with variance $1$.
We then compute the inverse variables $\lambda_i$ and $w_i$ of $T$
and test the algorithms with these inputs:
the norming constants $w_i$ typically span several orders of magnitude.
There are cases where either algorithm outperforms the other,
but in the average BI fares decisively better than BG.
In the worksheet, we repeated this
experiment $40$ times with dimension $n = 40$, working with $12$
significant digits; errors were measured as in the previous section.
The run is declared a failure if the error exceeds $0.1$:
there were two runs where both BG and BI failed,
another $30$ failed runs for BG and none other for BI.


We next consider matrices near $T_0$,
the Jacobi matrix with diagonal entries
equal to $0$ and off-diagonal entries equal to $1$:
for our purposes, $T_0$ is as good as the free Laplacian.
It turns out that $T_0$ is special from several points of view:
the values of both $\lambda_i$ and $w_i$ can be obtained explicitly,
there are no small gaps in the spectrum and all norming constants
$w_i$ have roughly the same size.
These features, particularly the last one, seem to favor BG.
Indeed, for $n = 40$, working with $12$ digits,
Gaussian perturbations of $T_0$ with small variance in
the off-diagonal entries favors BG:
among $40$ examples, there are no failures of either algorithm
but the errors are smaller for BG than for BI.


\vfill\eject

\bigskip\bigskip\bigbreak

{

\parindent=0pt
\parskip=0pt
\obeylines

Ricardo S. Leite, Departamento de Matemática, UFES 
Av. Fernando Ferrari, 514, Vitória, ES 29075-910, Brazil

Nicolau C. Saldanha and Carlos Tomei, Departamento de Matem\'atica, PUC-Rio
R. Marqu\^es de S. Vicente 225, Rio de Janeiro, RJ 22453-900, Brazil

\smallskip

rsleite@cce.ufes.br
nicolau@mat.puc-rio.br; http://www.mat.puc-rio.br/$\sim$nicolau/
tomei@mat.puc-rio.br

}


\begin{thebibliography}{[10]}

\bibitem{BG}{ de Boor, C. and Golub, G. H.,
{The numerically stable reconstruction of a Jacobi matrix
from spectral data},
Linear Algebra and its Applications 21, 245-260, 1978.}
\bibitem{BFR}{ Bloch, A. M., Flaschka, H. and Ratiu, T.,
{A convexity theorem for isospectral manifolds of Jacobi matrices
in a compact Lie algebra},
Duke Math. J., 61, 41-65, 1990.}
\bibitem{EP}{Erra, R. and Phillipe, B.,
{On some structured inverse eigenvalue problems},
Numerical Algorithms, 15, 15-35, 1997.}
\bibitem{GW}{ Gray, L. J. and Wilson, D. G.,
{Construction of Jacobi matrix from spectral data},
Linear Algebra Appl. 14, 131-134, 1976.}
\bibitem{Hald}{ Hald, O. H.,
{Inverse eigenvalue problems for Jacobi matrices},
Linear Algebra Appl. 14, 63-85, 1976.}
\bibitem{Hochstadt}{ Hochstadt, H.,
{On some inverse problems in matrix theory},
Arch. Math. 18, 201-207, 1967.}
\bibitem{Hochstadt2}{ Hochstadt, H.,
{On the construction of a Jacobi matrix from spectral data},
Linear Algebra Appl. 8, 435-446, 1974.}
\bibitem{LST}{ Leite, R. S., Saldanha, N. C. and Tomei, C.,
{New inverse data for tridiagonal matrices and
the asymptotics of Wilkinson's shift iteration},
preprint, {\tt www.arxiv.org/abs/math.NA/0412493}.}
\bibitem{Moser}{ Moser, J.,
{Finitely many mass points on the line under the influence
of an exponential potential},
In: Dynamic systems theory and applications, (ed. J.~Moser)
467-497, New York, 1975.}
\bibitem{Tomei}{ Tomei, C.,
{The Topology of Manifolds of Isospectral Tridiagonal Matrices},
Duke Math. J., 51, 981-996, 1984.}


\end{thebibliography}
\end{document}